\documentclass[11pt]{amsart}
\usepackage{amssymb,amsmath,amsthm}
\usepackage{a4}
\usepackage{epsf}
\usepackage{epsfig}
\begin{document}

\newcommand{\NP}{$\mathcal{N}\mathcal{P}$}
\newcommand{\newsymb}{{\mathcal P}}
\newcommand{\Pd}{{\mathcal P}^d}
\newcommand{\R}{{\mathbb R}}
\newcommand{\N}{{\mathbb N}}
\newcommand{\QQ}{{\mathbb Q}}
\newcommand{\ZZ}{{\mathbb Z}}
\newcommand{\STAB}{\mathrm{STAB}}

\newcommand{\enorm}[1]{\Vert #1\Vert}
\newcommand{\inter}{\matrm{int}}
\newcommand{\conv}{\mathrm{conv}}
\newcommand{\aff}{\mathrm{aff}}
\newcommand{\lin}{\mathrm{lin}}
\newcommand{\cone}{\mathrm{cone}}

\newcommand{\dist}{\mathrm{dist}}
\newcommand{\trans}{\intercal}
\newcommand{\diam}{\mathrm{diam}}
\newcommand{\pp}{\mathfrak{p}}
\newcommand{\pf}{\mathfrak{f}}
\newcommand{\pg}{\mathfrak{g}}
\newcommand{\PP}{\mathfrak{P}}
\newcommand{\pl}{\mathfrak{l}}
\newcommand{\pv}{\mathfrak{v}}
\newcommand{\cl}{\mathrm{cl}}
\newcommand{\bx}{\overline{x}}

\def\ip(#1,#2){#1\cdot#2}

\newtheorem{theorem}{Theorem}[section]
\newtheorem{corollary}[theorem]{Corollary}
\newtheorem{lemma}[theorem]{Lemma}
\newtheorem{remark}[theorem]{Remark}
\newtheorem{definition}[theorem]{Definition}  
\newtheorem{proposition}[theorem]{Proposition}  
\newtheorem{claim}[theorem]{Claim}
\numberwithin{equation}{section}

\title[Polynomial inequalities representing polyhedra]{Polynomial inequalities representing polyhedra${}^\star$} 

\thanks{$\star$ Supported by the DFG Research Center
        ``Mathematics for key technologies''
        (FZT 86) in Berlin.}
\author{Hartwig Bosse}
\address{Hartwig Bosse,Technische Universit\"at Berlin, Sekr. MA6-2, 
         Stra{\ss}e des 17.~Juni 135, D-10623 Berlin, Germany}
\email{bosse@math.tu-berlin.de}

\author{Martin Gr\"otschel}
\address{Martin Gr\"otschel, Konrad-Zuse-Zentrum f\"ur
  Informationstechnik (ZIB), 
         Taku\-str.~7, D-14195 Berlin-Dahlem, Germany}
\email{groetschel@zib.de}

\author{Martin Henk}
\address{Martin Henk, Universit\"at Magdeburg, Institut f\"ur Algebra und Geometrie,
  Universit\"ats\-platz 2, D-39106 Magdeburg, Germany} 
\email{henk@math.uni-magdeburg.de}


\begin{abstract}
Our main result is that every $n$-dimensional polytope can be described by at
most $(2n-1)$ polynomial inequalities and, moreover, these polynomials
can explicitly be constructed. For an $n$-dimensional pointed
polyhedral cone we prove the bound $2n-2$ and for arbitrary polyhedra
we get a constructible representation by $2n$ polynomial inequalities. 
\end{abstract}

\maketitle

\section{Introduction}
By a striking result of Br\"ocker and Scheiderer (see
\cite{Scheiderer:stability}, \cite{broecker:semialgebraic}, 
 \cite{bochnakcosteroy:real_algebraic_geometry}  and \cite{Mahe:broecker_scheiderer}), every basic closed semi-algebraic set  of the form      
\begin{equation*}
      \mathcal{S}=\left\{ x\in \R^n : \pf_1(x)\geq 0,\dots, \pf_l(x)\geq 0\right\}, 
\end{equation*}
where $\pf_i\in \R[x]$, $1\leq i\leq l$,  are
polynomials, can be represented by at most $n(n+1)/2$ polynomials,
i.e., there exist polynomials 
$\pp_1,\dots,\pp_{n(n+1)/2}\in\R[x]$ such that  
 \begin{equation*} 
     \mathcal{S}=
     \left\{ x\in \R^n : \pp_1(x)\geq 0,\dots, \pp_{n(n+1)/2}(x)\geq 0\right\}.
\end{equation*}
Moreover, in the case of basic open semi-algebraic sets, i.e., $\geq$ is
replaced by strict inequality, one can even bound the maximal number
of polynomials needed by the dimension $n$ instead of
$n(n+1)/2$. Rephrasing the results in terms of semi-algebraic
geometry, the {\em stability index} of every basic closed or open
semi-algebraic set is $n(n+1)/2$ or $n$, respectively. Both bounds are best possible. 

No explicit constructions, however, of such  
systems of polynomials are known nor whether the upper bound $n(n+1)/2$ can
be improved for semi-algebraic sets having additional structure such as 
convexity. Even in the very special case of
$n$-dimensional polyhedra almost nothing was known.   In \cite[Example
2.10]{broecker:semialgebraic} or in \cite[Example
4.7]{andradas_broecker_ruiz:constructible} 
a description of a regular convex $m$-gon in the plane 
by two polynomials is given. This result was generalised to arbitrary
convex polygons and three polynomial inequalities by vom Hofe
\cite{vom_Hofe:disser}.   Bernig \cite{Bernig:diplom} proved that, for $n=2$,  every  convex polygon
can  even be represented by two polynomial inequalities. In
\cite{groetschel_henk:poly_poly} a construction of $O(n^n)$ polynomial
inequalities  representing an $n$-dimensional simple
polytope is given. Based on ideas from \cite{Bosse:diplom}, 
here we give, in particular,  an explicit construction of 
$(2n-1)$ polynomials describing an arbitrary $n$-dimensional polytope.
Hence the general upper bound of $n(n+1)/2$ polynomials can be
improved (at least) for $n$-dimensional polytopes, and we
conjecture that the dimension itself is the right value  for this
special case.

In order to state our results we fix some notation. 
A polyhedron $P\subset \R^n$ is the intersection of 
finitely many closed halfspaces, i.e., we can write it as 
\begin{equation*}
  P=\left\{x\in\R^n : \ip(a_i,x) \leq b_i, \, 1\leq i\leq
    m\right\},
\end{equation*} 
for some $a_i\in\R^n$, $b_i\in\R$. Here $\ip(a,x)$ denotes the
standard inner product on $\R^n$. If $P$ is bounded then it is called
a polytope.  A pointed polyhedral cone $C\subset\R^n$ with apex at the
origin is the intersection of finitely many closed halfspaces of the
type 
\begin{equation*}
  C=\left\{x\in\R^d : \ip(a_i,x) \leq 0, \, 1\leq i\leq
    m\right\},
\end{equation*} 
$a_i\in\R^n$. For polynomials $\pp_i\in \R[x]$, $1\leq i\leq l$, we
denote by  
\begin{equation*}
 \newsymb(\pp_1,\dots,\pp_l) := \left\{ x\in \R^n : \pp_1(x)\geq
   0,\dots, \pp_{l}(x)\geq 0\right\} 
\end{equation*}
the associated basic closed semi-algebraic set generated by the polynomials.
\begin{theorem} Let $C\subset \R^n$ be an $n$-dimensional pointed
  polyhedral cone. Then we can construct $(2n-2)$ polynomials $\pp_i\in\R[x]$, $1\leq i \leq 2n-2$,
  such that $C=\newsymb(\pp_1,\dots,\pp_{2n-2})$.
\label{thm:cones}
\end{theorem}
The case of polytopes can be derived as  a consequence of the
construction behind Theorem \ref{thm:cones} and here we get 
\begin{theorem} Let $P\subset \R^n$ be an $n$-dimensional
  polytope. Then we can construct $(2n-1)$ polynomials $\pp_i\in\R[x]$, $1\leq i \leq 2n-1$,
  such that $P=\newsymb(\pp_1,\dots,\pp_{2n-1})$.
\label{thm:polytopes}
\end{theorem}
At the end of  Section 3  (see Definition \ref{def:polynomials}) 
we will give an explicit description of the
polynomials we employ. 
The construction behind the proof of Theorem
\ref{thm:polytopes} or Theorem \ref{thm:cones} can also be applied to
the interior of a polytope or a cone which
are open semi-algebraic sets. Furthermore, in  \cite[Proposition
2.5]{groetschel_henk:poly_poly} it is shown how a representation of a
polytope by polynomial inequalities can be used  to get a
representation of a polyhedron by polynomials. Applying this
proposition to Theorem \ref{thm:polytopes}  leads to 
\begin{corollary} Let $P\subset \R^n$ be an $n$-dimensional
  polyhedron. Then we can construct $2n$ polynomials $\pp_i\in\R[x]$, $1\leq i \leq 2n$,
  such that $P=\newsymb(\pp_1,\dots,\pp_{2n})$.
\end{corollary}

The paper is organised as follows. In Section 2 we give, for a pointed
cone $C$, a construction  of two polynomials $\pp_{C,\varepsilon}, \pp_0$
such that $C$ is ``nicely  approximated'' by
$\newsymb(\pp_{C,\varepsilon}, \pp_0)$. Then, for a face $F=C\cap
\{x\in\R^n : \ip(a_i,x)=0,\, i\in I_F\}$ of $C$,  
we apply this construction to the cone $C_F=\{x\in\R^n : \ip(a_i,x)
\leq 0, i\in I_F\}$, where $I_F$ denotes the index set of active
constraints of $F$. In that way we get an approximation of $C_F$ by a semi-algebraic set of the type $\newsymb(\pp_{C_F,\varepsilon}, \pp_F)$.  In Section 3 we study the relations between the set $\newsymb(\pp_{C_{F\cap G},\varepsilon}, \pp_{F\cap G})$ and 
$\newsymb(\pp_{C_F,\varepsilon}, \pp_F)$, $\newsymb(\pp_{C_G,\varepsilon}, \pp_G)$ for two different faces $F$ and $G$ of the same dimension.  
Thereby, it turns out that we may multiply all polynomials $\pp_{C_F,\varepsilon}$ belonging to faces of the same dimension as well as the   polynomials $\pp_{F}$ in order
 to get a representation of a
pointed polyhedral cone by polynomials. In Section 4 we give a brief
outlook why we are
interested in such a polynomial representation of polytopes and what
might be achievable by such a representation with respect to hard combinatorial optimisation problems. 

\section{Approximating cones}
In the following we use some standard terminology and facts from the
theory of polyhedra for which we refer to the books \cite{McMShe:upper} and 
\cite{ziegler:polytopes}.
For the  approximation of a cone by a closed semi-algebraic
set consisting of two polynomials we need a lemma about the
approximation of a polytope by a strictly convex polynomial which was
already shown in \cite[Lemma 2.6]{groetschel_henk:poly_poly}. Since it
is essential for the explicit construction of the polynomials we state it
here. To this end, let $B^n$ be  the $n$-dimensional unit ball
centred at the origin. The diameter of a polytope is denoted by
$\diam(P)$, i.e., $\diam(P)=\max\{\enorm{x-y}: x,y\in P\}$, where
$\enorm{\cdot}$ denotes the Euclidean norm. 
\begin{lemma} Let $P=\{x\in\R^n : \ip(a_i,x) \leq b_i,\,\, 1\leq i\leq
  m\}$ be an $n$-dimensional polytope. For $1\leq i\leq m$ let 
$$ 
  \pv_i(x):=
   \frac{2\ip(a^i,x) - h(a_i)+h(-a_i)}
        {h(a_i)+h(-a_i)}, 
$$
where $h(a):=\max\{\ip(a,x) : x\in P\}$ is the support function of
$P$. Let $\varepsilon>0$, choose an integer $k$ such that $k >
\ln(m)/(2\ln(1+\frac{2\varepsilon}{(n+1)\diam(P)}))$, and set  
$$
  \pp_{P,\varepsilon}(x):= \sum_{i=1}^m
  \frac{1}{m}\,\left[\pv_i(x)\right]^{2\,k}\quad\text{ and }\quad 
  K_\varepsilon:=\{x\in\R^n : \pp_{P,\varepsilon}(x)< 1\}.
$$
Then we have 
$P\subset K_\varepsilon\subset P+\varepsilon\,B^n$.
\label{lem:approximating} 
\end{lemma} 
\begin{proof} \cite[Lemma 2.6]{groetschel_henk:poly_poly}.
\end{proof}

 Now let 
\begin{equation} 
  C=\left\{x\in\R^n : \ip(a_i,x) \leq 0, \, 1\leq i\leq
    m\right\},
\label{eq:cone}
\end{equation} 
be a pointed $n$-dimensional cone with $\enorm{a_i}=1$, $1\leq i\leq m$. The set of all
$k$-dimensional faces ($k$-faces for short) is denoted by
$\mathcal{F}_k$, $0\leq k\leq n-1$. For a $k$-face $F$, we denote by 
$I_F:=\{i: \ip(a_i,x)=0\,\text{ for all } x\in F\}$ the set of active
constraints.
We always assume that our 
representation \eqref{eq:cone} of $C$ is non-redundant, hence $\{x\in C : \ip(a_i,x) =
0\}$ is an $(n-1)$-face (facet) of $C$ for $1\leq i\leq m$.  For each $F$, let 
\begin{equation}
   u_F:=\frac{\sum_{i\in I_F} a_i}{\Vert\sum_{i\in I_F} a_i\Vert}\quad
   \text{and}\quad \pp_F(x):=-\ip(u_F,x).
\label{eq:simple_polynomials}
\end{equation}
$u_F$ is an outer unit normal vector of the face $F$, i.e.,
$F=C\cap\{x\in\R^n : \pp_F(x)=0\}$ and $C\setminus F\subset
\{x\in\R^n : \pp_F(x)>0\}$. The only vertex, i.e., $0$-face, of $C$ is the
origin, and in this case, we denote the above  outer unit normal vector
and the polynomial 
by $u_0$ and  $\pp_0$, respectively. 
In the next lemma we construct a basic closed semi-algebraic
set consisting of two polynomials that gives  a nice and
controllable approximation of $C$.  In what follows we will often use
some constants depending on  the cone or
polytope. All of these constants are  explicitly  computable by
elementary methods, but in order to keep the presentation simple we do not go into the
details here.

\begin{lemma} For every $\varepsilon\in(0,1/2]$ we can construct a polynomial
  $\pp_{C,\varepsilon}(x)$ such that 
\begin{equation*}
\begin{split}
{\rm i)}& \,\,\left\{x+\varepsilon\,(\ip(u_0,x)) B^n : x\in
  C\right\} \subset \newsymb(\pp_{C,\varepsilon},\pp_0)\subset
\left\{x+\omega_C\,\varepsilon\,(\ip(u_0,x)) B^n : x\in C\right\},\\
 {\rm ii)}& \,\,\left\{x\in\R^n : \pp_{C,\varepsilon}(x) = 0, \pp_0(x) = 0\right\}
 =\{0\}, \\
 {\rm iii)}&\,\, \left\{x+\varepsilon\,(\ip(u_0,x)) B^n : x\in
  C,\,\pp_0(x) >0 \right\} \subset \left\{x\in\R^n :
  \pp_{C,\varepsilon}(x)> 0\right\}, 
\end{split}
\end{equation*}
where $\omega_C \geq 1$ is a constant depending only on $C$.
\label{lem:cone_approximation}
\end{lemma}


\begin{proof}
  Firstly, observe that for $n=1$ there is nothing to do, because we
  may set $\pp_{C,\varepsilon}(x):=\pp_0(x)$ and $\omega_C=1$, say. So let
  $n\geq 2$. For ease of notation we may assume that    
  $-u_0=e_n$, the $n$-th unit vector, which can be achieved by a
  suitable rotation. Due to this choice $C\cap\{x\in\R^n: x_n=1\}$ is
  an $(n-1)$-dimensional polytope $P$, which  we identify  with its
  image under the orthogonal projection onto
  $\R^{n-1}$. Thus let 
$P=\{x\in\R^{n-1} : \ip(\widetilde{a}_i,x)\leq
  \widetilde{b}_i,\, 1\leq i\leq m\}$, for some $\widetilde{a}_i\in\R^{n-1}$,  
  $\enorm{\widetilde{a}_i}=1$, $\widetilde{b}_i\in\R$. With this
  notation  we may write $C$ as the
  homogenisation of $P$, i.e., 
  $C=\{x_n\,(x,1)^\intercal : x\in P,\,x_n\geq 0\}$. For $\mu\geq 0$
  let 
\begin{equation*}
  P_\mu = \{x\in\R^{n-1} : \ip(\widetilde{a}_i,x)\leq
  \widetilde{b}_i+\mu,\, 1\leq i\leq m\}.
\end{equation*}    
Then 
\begin{equation*}
 P+\mu\,B^{n-1}\subset P_\mu\subset P+\omega_P\,\mu\,B^{n-1},
\end{equation*}
for a certain constant $\omega_P\geq 1$ depending only on
$P$. From Lemma \ref{lem:approximating}  we get that, for
every $\nu>0$, we can construct  a strictly convex polynomial $\pp_{P_{\mu,\nu}}$
such that 
\begin{equation}
   P_\mu\subset \left\{x\in\R^{n-1} :
     \pp_{P_{\mu,\nu}}(x)<1\right\}\subset P_\mu+\nu\,B^{n-1}.
\label{eq:poly_app}
\end{equation}
In particular, $\pp_{P_{\mu,\nu}}$ can be written as 
$
  \pp_{P_{\mu,\nu}}(x)=\sum_{i=1}^m \lambda_i[\ip(\widetilde{a}_i,x)-\alpha_i]^{2k}$
for certain constants $\lambda_i\in\R_{>0}$, $\alpha_i\in\R$,
$k\in\N$, depending on $P_\mu$ and $\nu$ (cf.~Lemma \ref{lem:approximating}). 
For a scalar  $x_n > 0$ we immediately get 
\begin{equation}
\begin{split}
  x_n\,P_\mu&\subset \{x\in\R^{n-1}: \sum_{i=1}^m
  \lambda_i[\ip(\widetilde{a}_i,x)-x_n\alpha_i]^{2k}<(x_n)^{2k}\}\\&\subset x_n\,P+x_n(\nu+\omega_P\mu)\,B^{n-1}.
\end{split}
\label{eq:inclus}
\end{equation}
Since $\widetilde{a}_1,\ldots,\widetilde{a}_m$ are the outer normal vectors 
of an $(n-1)$-dimensional polytope, these inclusions hold for $x_n=0$ as well, 
if we replace $<$ by $\leq$. Hence, with  
\begin{equation*}
\overline{\pp}_{P_{\mu,\nu}}(x)=(x_n)^{2k}-\sum_{i=1}^m
\lambda_i[\ip(\widetilde{a}_i,{(x_1,\dots,x_{n-1})}^\intercal)-x_n\alpha_i]^{2k} 
\end{equation*}
and $\pp_0(x)=x_n$, for $x=(x_1,\dots,x_n)^\intercal\in\R^n$,   we get 
\begin{equation}
 \begin{split}
  {\rm i)} \quad & \{x\in\R^n : \overline{\pp}_{P_{\mu,\nu}}(x)=0, \, 
  \pp_0(x)=0\} = \{0\},\\
  {\rm ii)} \quad & x_n\,P_\mu \subset \{x\in\R^n : \overline{\pp}_{P_{\mu,\nu}}(x)\geq 0\}, \text{ for } x_n\geq 0,\\
  {\rm iii)} \quad & x_n\,P_\mu \subset \{x\in\R^n : \overline{\pp}_{P_{\mu,\nu}}(x)>0\}, \text{ for } x_n>0.
\end{split}
\label{eq:prop}
\end{equation}
From \eqref{eq:inclus} we conclude that 
\begin{equation}
   \newsymb(\overline{\pp}_{P_{\mu,\nu}},\pp_0)
   \subset\left\{x+x_n(\nu+\omega_P\mu)B^n,\,x\in C\right\}. 
\label{eq:good_app}
\end{equation}
With  $\gamma=\max\{(1-\ip(a_i,e_n))^{-1/2}: 1\leq i\leq m\}$ and
by some elementary calculations we get    for 
$y\in \left\{x+x_n\left(\frac{\mu}{\mu+\gamma}\right)B^n,\,x\in
     C\right\}$ that 
\begin{equation}
   (y_1,\dots,y_{n-1})^\intercal \in  y_n P_\mu.
\label{eq:lower_incl}
\end{equation}
Thus we have by \eqref{eq:prop} ii) 
 \begin{equation}
 \left\{x+x_n\left(\frac{\mu}{\mu+\gamma}\right)B^n,\,x\in
     C\right\}\subset  \newsymb(\overline{\pp}_{P_{\mu,\nu}},\pp_0).
\label{eq:not_too_good}
\end{equation}
Now, for a given $\varepsilon\in(0,1/2]$,  we may choose $\mu$ and $\nu$ such that 
$\mu/(\mu+\gamma)=\varepsilon$
and $\nu+\omega_P\mu\leq 4\gamma\omega_P\varepsilon$. With $\omega_C:=4\gamma\omega_P$ and 
 $\pp_{C,\varepsilon}:=\overline{\pp}_{P_{\mu,\nu}}$
for this special choice
of parameters we get by \eqref{eq:good_app} and
\eqref{eq:not_too_good} the statement i) of the lemma. Property 
ii) is an immediate consequences of \eqref{eq:prop} i) and the last statement follows from \eqref{eq:lower_incl} and \eqref{eq:prop} iii).
\end{proof}

\begin{remark} \hfill
\begin{enumerate}
\item[{\rm i)}] 
The main geometric message of Lemma  \ref{lem:cone_approximation} is that we can construct a cone of the  type 
  $\newsymb(\pp_{C,\varepsilon},\pp_0)$, which  is not too far away from $C$,
but at the same time we also know that
$\newsymb(\pp_{C,\varepsilon},\pp_0)$ is not too close to $C$. This 
property of $\newsymb(\pp_{C,\varepsilon},\pp_0)$ plays a key role in our construction. 
\item[{\rm ii)}]
As constant $\omega_P$ in the above proof we can take $R(P)/r(P)$, where $R(P)$ and $r(P)$ denote the radii of two concentric balls such that   $x+r(P)\,B^n\subset P\subset x+R(P)\,B^n$.
\end{enumerate}
\end{remark}
For a $k$-face $F$ of $C$, let $C_F=\{x\in\R^n : \ip(a_i,x)\leq
0,\, i\in I_F\}$ be the face-cone of $F$.  $C_F$ is an $n$-dimensional cone containing a $k$
dimensional linear subspace, namely $\lin(F)$, the linear hull of $F$.
The $(n-k)$-dimensional orthogonal complement $\lin(F)^\perp$ of
$\lin(F)$ is given by $\lin\{a_i : i\in I_F\}$. If we apply the
construction of Lemma \ref{lem:cone_approximation} to $C_F\cap
\lin(F)^\perp$ (in the space $\lin(F)^\perp$) we get a generalisation
of  Lemma \ref{lem:cone_approximation} from the face-cone of the vertex to arbitrary $k$-faces of $C$. 
\begin{corollary} Let $F$ be a $k$-face of $C$ with $0\leq k\leq
  n-1$. For every $\varepsilon\in(0,1/2]$  we can construct a polynomial
  $\pp_{C_F,\varepsilon}(x)$ such that   
\begin{equation*}
\begin{split}
{\rm i)}& \,\,\left\{x+\varepsilon\,(\ip(u_F,x)) B^n : x\in
  C_F\right\} \subset \newsymb(\pp_{C_F,\varepsilon},\pp_F) \\&\,\,\,
 \hphantom{ \left\{x+\varepsilon\,(\ip(u_F,x)) B^n : x\in
  C_F\right\}}\subset
\left\{x+\omega_{C_F}\,\varepsilon\,(\ip(u_F,x)) B^n : x\in C_F\right\},\\
 {\rm ii)}& \,\,\left\{x\in\R^n : \pp_{C_F,\varepsilon}(x) = 0, \pp_F(x) = 0\right\}
 =\lin(F), \\
 {\rm iii)}&\,\, \left\{x+\varepsilon\,(\ip(u_F,x)) B^n : x\in
  C_F,\,\pp_F(x) >0 \right\} \subset \left\{x\in\R^n :
  \pp_{C_F,\varepsilon}(x)> 0\right\},
\end{split}
\end{equation*}
where $\omega_{C_F}\geq 1$ is a constant depending only on $C$.
\label{cor:cone_approximation_general}
\end{corollary}
We note that, for a facet $F$ of $C$ and $\varepsilon\in (0,1/2]$, we just have 
 (cf. proof of Lemma \ref{lem:cone_approximation}) 
\begin{equation}
\pp_{C_F,\varepsilon}(x)=\pp_F(x)=\ip(-u_F,x).
\label{eq:facets}
\end{equation}

\section{Multiplying polynomial inequalities} 
The main objective of our proof strategy is to multiply, for each $k\in\{0,\dots,n-1\}$, all the polynomials
$\pp_{C_F,\varepsilon}$, $F\in\mathcal{F}_K$, and $\pp_F$,
$F\in\mathcal{F}_k$, such that for a special choice  of the parameters
$\varepsilon$, the
arising $2n$ polynomials give a complete description of the cone $C$. To this
end, we have to study, for two $k$-faces $F$ and $G$,
the relations between $\newsymb(\pp_{C_F,\varepsilon},\pp_F)$,
$\newsymb(\pp_{C_G,\varepsilon},\pp_G)$, and $\newsymb(\pp_{C_{F\cap
    G},\varepsilon},\pp_{F\cap G})$.

\begin{lemma} Let $F,G$ be $k$-faces of $C$ and let
  $\varepsilon_k\in(0,1/2]$. Then we can find an  $\varepsilon_{F,G}\in (0,1/2]$ such
  that 
\begin{align*}
  \big\{x+\varepsilon_{F,G}\,(\ip(u_{F\cap G},x)) B^n &: x\in
  C_{F\cap G}\big\}  \subset \\ 
  &\left\{x+\varepsilon_k\,(\ip(u_F,x))\, B^n : x\in
  C_{F},\,\ip(-u_F,x) >0\right\} \\
  \cup &\left\{x+\varepsilon_k\,(\ip(u_G,x))\, B^n : x\in
  C_{G},\, \ip(-u_G,x)>0\right\} \\
   \cup &\left(\lin(F)\cap\lin(G)\right).
\end{align*}
\label{lem:inclusion}
\end{lemma}
\begin{proof} Let $C_{F\cap G} = \lin(F\cap G)
+ \cone\{v_1,\dots,v_r\}$ for some points $v_i\in \lin(F\cap G)^\perp$,
where $\cone$ denotes the conical  hull. Since both, 
$\frac{1}{2}(u_F+u_G)$ and $u_{F\cap G}$, are  outer normal vectors of the face $F\cap G$ we find that
\begin{equation*} 
\rho=\min\left\{ \frac{\ip(\frac{1}{2}(u_F+u_G),v_i) }{
    \ip(u_{F\cap G},v_i)}: 1\leq i \leq r \right\} >0.
\end{equation*}
Hence, for $x\in C_{F\cap G}$, we get  
\begin{equation} 
  \max\{\ip(-u_F,x), \ip(-u_G,x)\} \geq
  \ip(\frac{1}{2}(-u_F-u_G),x) \geq \rho  \ip((-u_{F\cap G}),x).
\label{eq:max}
\end{equation}
If $\ip(u_{F\cap G},x)=0$ then $x\in\lin(F\cap G)\subset \lin(F)\cap\lin(G)$.
Otherwise we have $\ip(-u_{F\cap G},x)>0$, and with $\varepsilon_{F,G}:=\min\{ \rho\varepsilon_k,1/2\}$
and \eqref{eq:max} we get the required inclusion. 

\end{proof}
 
As a corollary we get that we can find $\varepsilon_k$, $0\leq k\leq n-1$, such
  that a  cone of the type $\newsymb(\pp_{C_{F\cap G}, \varepsilon_{\dim(F\cap
      G)}},\,\pp_{F\cap G})$, $F,G\in\mathcal{F}_k$, is covered by the interior of
  $\newsymb(\pp_{C_F,\varepsilon_k}, \pp_F)$, the interior of
  $\newsymb(\pp_{C_G,\varepsilon_k}, \pp_G)$, and the linear space
  $\lin(F)\cap\lin(G)$.
\begin{corollary} We can determine positive constants $\varepsilon_k\leq 1/2$, $0\leq k\leq n-1$, such that 
for any pair of two different $k$-faces $F$ and $G$ of $C$, $k\in\{0,\dots,n-1\}$, 
\begin{equation}
\begin{split}
   \mathcal{P}(&\pp_{C_{F\cap G}, \varepsilon_{\dim(F\cap G)}},\,\pp_{F\cap
     G})  \subset \\ &\left\{x\in\R^n : \pp_{C_F,\varepsilon_k}(x)> 0, \pp_F(x)>0 \right\} \\
    \cup  &\left\{x\in\R^n : \pp_{C_G,\varepsilon_k}(x)> 0, \pp_G(x)>0 \right\} \\
    \cup &\left\{x\in\R^n :  \pp_{C_F,\varepsilon_k}(x)= 0, \pp_F(x)=0, \pp_{C_G,\varepsilon_k}(x)= 0, \pp_G(x)=0 \right\}.
\end{split}
\end{equation}
\label{cor:intersection}
\end{corollary}
\begin{proof}  By \eqref{eq:facets} we may set  $\varepsilon_{n-1}:=1/2$ and in view of  Corollary \ref{cor:cone_approximation_general} and
  Lemma \ref{lem:inclusion} we just have to say how to calculate the 
  numbers $\varepsilon_k$, $0\leq k\leq n-2$.  For two faces $F,G \in\mathcal{F}_k$ 
  the 
  proof of Lemma \ref{lem:inclusion} (the $\varepsilon_{F,G}$ constructed
  there) leads to an upper bound on
  $\varepsilon_{\dim(F\cap G)}$ provided we know $\varepsilon_k$. Hence, for $k=n-2,\dots,0$, we can calculate suitable numbers $\varepsilon_k$ via  
\begin{equation*}
  \varepsilon_k := \min_{k+1\leq l\leq n-1}\min_{ F,G\in \mathcal{F}_l}\left\{\varepsilon_{F,G}:
    \dim(F\cap G)=k \right\}.
\end{equation*}
 \end{proof}

Since every $(n-2)$-face $H$ of $C$ is given 
  by the intersection of two uniquely determined facets $F$ and $G$ of
  $C$ we may even set (cf.~\eqref{eq:facets}) 
\begin{equation}
   \varepsilon_{n-2}:= 1/2,\quad \pp_{C_H,\varepsilon_{n-2}}(x) := \pp_H(x)=\ip(-u_H,x)
\label{eq:ridges}
\end{equation}
without violating the validity of Corollary \ref{cor:intersection}. 

Now we come to the definition of the polynomials, which give us a
representation of an $n$-dimensional pointed polyhedral cone and to the proofs of Theorem \ref{thm:cones} and Theorem \ref{thm:polytopes}.

\begin{definition} Let $\varepsilon_k$, $0\leq k\leq n-1$, be chosen
  according to Corollary \ref{cor:intersection} and \eqref{eq:ridges}. For $F\in\mathcal{F}_k$, let $\pp_F, \pp_{{C_F},\varepsilon_k}\in\R[x]$ be given as in \eqref{eq:simple_polynomials}, Lemma \ref{lem:cone_approximation}, \eqref{eq:facets}, and \eqref{eq:ridges}. Then, 
 for $k=0,\dots,n-1$, let 
\begin{equation*}
 \PP_{k,1}(x) := \prod_{F\in\mathcal{F}_{k}} \pp_F(x)  \quad\text{ and }\quad \PP_{k,2}(x) := \prod_{F\in\mathcal{F}_{k}} \pp_{C_F,\varepsilon_k}(x).
\end{equation*}
\label{def:polynomials}
\end{definition}

\begin{proof}[Proof of Theorem \ref{thm:cones}] First we show  that 
\begin{align*}
  C=\big\{x\in\R^n : \PP_{k,1}(x)\geq 0, \,\PP_{k,2}(x)\geq 0,\, k=0,\dots,n-1 \big\}.
\end{align*}
The inclusion $\subset$ is obvious. So let $y\notin C$,
  but suppose that $y$ satisfies all the polynomial
  inequalities. Since $y\notin C$ one of the facet defining
  inequalities has to be violated, i.e., there exists an $(n-1)$-face
  $F$ with $\pp_F(y)<0$. Hence we may define
  $p\in\{0,\dots,n-1\}$ 
  as the minimum
  number (index) for which one of the factors in the polynomials
  $\PP_{p,1}(x)$ or $\PP_{p,2}(x)$ is violated. Since both, $\PP_{0,1}(x)$
  and $\PP_{0,2}(x)$, consists only of one polynomial we have
  $p\in\{1,\dots,n-1\}$. 

  Let $F\in \mathcal{F}_p$ such that $\pp_F(y)<0$ or
  $\pp_{C_F,\varepsilon_p}(y)<0$. Since $\PP_{p,1}(y)\geq 0$ and
  $\PP_{p,2}(y)\geq 0$ there must exist a $G\in \mathcal{F}_p$ with
  $\pp_G(y)\leq 0$ (in the case that $\pp_F(y)<0$)  or with
  $\pp_{C_G,\varepsilon_p}(y)\leq 0$ (if
  $\pp_{C_F,\varepsilon_p}(y)<0$). Thus we know that $y$ is neither 
  contained in the interior of the cone
  $\newsymb(\pp_{C_F,\varepsilon_{l}},\pp_F)$ nor in the interior of
  $\newsymb(\pp_{C_G,\varepsilon_{l}},\pp_G)$ nor in the linear space
  $\lin(F)\cap\lin(G)$.  By the choice of
  $\varepsilon_{\dim(F\cap G)}$ and Corollary  \ref{cor:intersection}, however, 
  those points $y$ are cut off by
  the cone  
  $\newsymb(\pp_{C_{F\cap G}, \varepsilon_{\dim(F\cap G)}},\,\pp_{F\cap G})$. Thus we must have   
\begin{equation*}
   y\notin \newsymb(\pp_{C_{(F\cap G)},\varepsilon_{\dim(F\cap
       G)}},\pp_{F\cap G})
\end{equation*}
contradicting the minimum property  of $p$. Finally, we observe that 
by \eqref{eq:facets} $\PP_{n-1,1}=\PP_{n-1,2}$, by
\eqref{eq:ridges} $\PP_{n-2,1}=\PP_{n-2,2}$ and hence we only have
$2n-2$ polynomials. 
\end{proof}

The  key to this algebraic proof are the special geometric properties 
{\rm i)} to {\rm iii)} 
of the approximative sets introduced in Corollary  
\ref{cor:cone_approximation_general}. 
These relations in combination with the result of Corollary \ref{cor:intersection} ensure that, for each pair of faces $F$, $G$, 
the set $\mathcal{P}(\pp_{C_{F\cap G}, \varepsilon_{\dim(F\cap G)}},\,\pp_{F\cap G})$ 
is contained in a special way in the union of the corresponding sets constructed for $F,G$ respectively, and this inclusion allows us to multiply those polynomials 
the latter are based on.


\begin{proof}[Proof of Theorem \ref{thm:polytopes}] Let $P\subset\R^n$
  be an $n$-dimensional polytope and let $C\subset\R^{n+1}$ be the
  $(n+1)$-dimensional 
  pointed polyhedral cone $C=\{x_{n+1}(x,1)^\intercal : x\in
  P\}$. Theorem \ref{thm:cones} shows that we construct $2n$
  polynomials describing $C$, where, in particular,  one polynomial
  ($\PP_{0,1}(x)$ in the notation of Definition \ref{def:polynomials})
  describes just a supporting hyperplane of $C$ at the origin.  
  Fixing the last coordinate to
  $x_{n+1}=1$ in these polynomials  gives a representation of $P$ by $2n$ polynomials. The
  polynomial $\PP_{0,1}(x)$, however, is apparently redundant for the polytope. 
\end{proof}

\begin{remark} We want to remark that for a
  polytope $P=\{x\in\R^n : \ip(a_i,x)\leq b_i,1\leq i\leq m\}$ with
  rational input data $a_i,\,b_i$ all the constants  involved 
  in the
  construction of the polynomials $\pp_{C_F,\varepsilon}$ can be
  substituted by certain rational numbers. Moreover, these numbers 
  can be calculated by well known methods from Linear Programming or
  Computational Geometry (cf.~\cite{Bosse:diplom}).  
\end{remark}

\section{Outlook}
The usual method 
to attack hard combinatorial optimisation problems is the polyhedral
approach. The basic idea here is a ``change of the representation'' of the problem, namely, to represent
combinatorial objects (such as the tours of a travelling salesman,
the independent sets of a matroid, or the stable sets in a graph)
as the vertices of a polytope. If one can find complete or
tight partial representations of polytopes of this type by linear
equations and inequalities,
linear programming (LP) techniques can 
be employed to solve the associated combinatorial optimisation
problem, see \cite{GLS:ellipsoid}. Even in the case where only
partial inequalities of the polyhedra associated with
  combinatorial
problems are known, LP techniques (such as cutting planes and column
generation) have resulted in very successful exact or approximate
solution methods. One prime example for this methodology is the
travelling salesman problem, see \cite{applebixbyetal:travel} and the corresponding
web page at {\tt http://www.math.princeton.edu/tsp/}.
Progress of the type  may also be possible via a 
``polynomial-representation approach''. Of course, since the degree of the 
polynomials in a such a polynomial representation is in general very
high (see e.g.~\cite[Proposition 2.1]{groetschel_henk:poly_poly}), and since polynomial inequalities are much
harder to treat than linear inequalities, we can not expect that such
an exact polynomial representation yields immediately a new method for combinatorial
optimisation problems. However, if we can answer questions like  how
well can we construct a small number of 
``simple'' polynomials $\pp_1,\dots,\pp_k$ 
such that a given  polytope  (or a general closed semi-algebraic set) is well approximated by the corresponding polynomials,
or how well can it be described or approximated by polynomials of
total degree $k$, then we believe that those results lead to 
a new approach to combinatorial optimisation 
problems via non-linear methods. 
We do know, of course, that these indications of possible future
results are mere speculation. Visions of this type, however,
were the starting point of the results presented in this
paper. 

\providecommand{\bysame}{\leavevmode\hbox to3em{\hrulefill}\thinspace}
\providecommand{\MR}{\relax\ifhmode\unskip\space\fi MR }
\providecommand{\MRhref}[2]{%
  \href{http://www.ams.org/mathscinet-getitem?mr=#1}{#2}
}
\providecommand{\href}[2]{#2}

\end{document}